\documentclass[12pt]{amsart}
\usepackage{amsmath,amscd,amssymb,amsfonts}
\def\ver{May 19, 2008, v.5}
\def\scirc{\raise.2ex\hbox{${\scriptstyle\circ}$}}
\def\csbull{\raise.4ex\hbox{${\scriptscriptstyle\bullet}$}}
\def\mopl{\hbox{$\bigoplus$}}
\def\msum{\hbox{$\sum$}}
\def\mcup{\hbox{$\bigcup$}}
\def\mcap{\hbox{$\bigcap$}}
\def\C{{\mathbb C}}
\def\D{{\mathbb D}}
\def\G{{\mathbb G}}

\def\P{{\mathbb P}}
\def\Q{{\mathbb Q}}
\def\Z{{\mathbb Z}}
\def\DD{{\mathcal D}}
\def\M{{\mathcal M}}
\def\O{{\mathcal O}}
\def\PP{{\mathcal P}}
\def\ZZ{{\mathcal Z}}
\def\ok{\overline{k}}
\def\oQ{\overline{\mathbb Q}}
\def\tD{\widetilde{D}}
\def\tX{\widetilde{X}}
\def\tcZ{\widetilde{\mathcal Z}}
\def\CH{\hbox{{\rm CH}}}
\def\Spec{\text{{\rm Spec}}\,}
\def\Im{\hbox{{\rm Im}}}
\def\Ker{\hbox{{\rm Ker}}}
\def\Hom{\hbox{{\rm Hom}}}
\def\Ext{\hbox{{\rm Ext}}}
\def\Hdg{\hbox{{\rm Hdg}}}
\def\MHM{\text{{\rm MHM}}}
\def\MHS{\text{{\rm MHS}}}
\def\Sym{\hbox{{\rm Sym}}}
\def\Gal{\text{{\rm Gal}}}
\def\Res{\hbox{{\rm Res}}}
\def\Aut{\hbox{{\rm Aut}}}
\def\Jac{\hbox{{\rm Jac}}}
\def\Pic{\hbox{{\rm Pic}}}
\def\Proj{\hbox{{\rm Proj}}}
\def\div{\hbox{{\rm div}}\,}
\def\Sing{\hbox{{\rm Sing}}\,}
\def\supp{\hbox{{\rm supp}}\,}
\def\codim{\hbox{{\rm codim}}\,}
\def\Gr{\hbox{{\rm Gr}}}
\def\ind{\text{\rm ind}}
\def\BM{{\rm BM}}
\def\DR{{\rm DR}}
\def\({{\rm (}}
\def\){{\rm )}}

\begin{document}
\title{Hodge-Type Conjecture for Higher Chow Groups}
\author{Morihiko Saito}
\address{RIMS Kyoto University, Kyoto 606-8502 Japan}
\email{msaito@kurims.kyoto-u.ac.jp}
\date{\ver}
\begin{abstract}
Let $X$ be a smooth quasi-projective variety over the algebraic
closure of the rational number field.
We show that the cycle map of the higher Chow group to Deligne
cohomology is injective and the higher Hodge cycles are generated by
the image of the cycle map as conjectured by Beilinson and Jannsen,
if the cycle map to Deligne cohomology is injective and the Hodge
conjecture is true for certain smooth projective varieties over the
algebraic closure of the rational number field.
We also verify the conjecture on the surjectivity in some cases of
the complement of a union of general hypersurfaces in a smooth
projective variety.
\end{abstract}
\maketitle

\centerline{{\bf Introduction}}

\bigskip\noindent
Let $X$ be a smooth projective variety over a subfield $k$ of $\C$,
and $\CH^p(X,n)_{\Q}$ the higher Chow group with $\Q$-coefficients
[8].
We have a cycle map to the Deligne cohomology of $X_{\C} :=
X\otimes_{k}\C$ (see [2], [9], [14], [15], [16], [19], [31], etc.):
$$
\CH^p(X,n)_{\Q}\to H_{\DD}^{2p-n}(X_{\C},\Q(p)).
\leqno(0.1)
$$
We are interested in its injectivity when $k$ is a subfield of $\oQ$.
For $n=0$, this is conjectured by Bloch and Beilinson [2] (at least
for cycles algebraically equivalent to zero), see also [20].
Note that the injectivity of (0.1) for smooth quasi-projective
varieties would imply the injectivity of the refined cycle map in
[32] (see also [1]).
It is expected that (0.1) would be bijective if we replace the
Deligne cohomology by a certain extension group in the derived
category of a conjectural category of mixed motives [4], and that
higher extension groups $\Ext^i$ of mixed motives over a number
field would vanish for $i>1$, see [2], [7], [20], etc.
Since the Deligne cohomology is expressed as an extension group in
the derived category of mixed Hodge structures, the problem is
closely related to the full faithfulness of the forgetful functor
from the (conjectural) category of mixed motives to that of mixed
Hodge structures.
Note that the latter problem may be viewed as an extension of the
Hodge conjecture which predicts the full faithfulness for pure
motives, see [12].

Let now $X$ be a smooth quasi-projective variety over $\oQ$.
Then the cycle map induces
$$
\CH^p(X,n)_{\Q}\to\Hom_{\MHS}(\Q, H^{2p-n} (X_{\C},\Q)(p)),
\leqno(0.2)
$$
where the target is called the group of higher Hodge cycles.
It is conjectured by Beilinson [3] and Jannsen [20] that (0.2)
would be surjective.
Here the source can be replaced with $\CH^p(X_{\C},n)_{\Q}$
by spreading cycles out as in Remark~(1.5)(ii).
Jannsen showed that the injectivity of (0.1) for $n=0$ is
essentially equivalent to the surjectivity of (0.2) for $n=1$,
where $X$ in (0.2) is the complement of the support of an algebraic
cycle of codimension $p$ on $X$ in (0.1).
However, any philosophical reason for the surjectivity of (0.2) does
not seem to have been known in general.
We show in this paper the following (see (4.1-2) for a more precise
statement):

\medskip\noindent
{\bf 0.3.~Theorem.} {\it If {\rm (0.1)} is injective and the Hodge
conjecture is true for any smooth projective varieties over $\oQ$,
then {\rm (0.2)} is surjective and {\rm (0.1)} is injective for any
smooth quasi-projective varieties over $\oQ$.
In particular, the refined cycle map {\rm [32] (}see also {\rm [1])}
is injective in this case.
}

\medskip
We prove this in a more general situation including the case
of systems of realizations [12], [13], [20].
Theorem~(0.3) gives evidence for Bloch's conjecture [6] for surfaces
with $p_{g}=0$, Murre's conjecture [25] on the Chow-K\"unneth
decomposition and Voisin's conjecture [35] on the countability of
$\CH_{\ind}^{2}(X,1)_{\Q}$.
Indeed, these can be reduced to the injectivity of the refined cycle
map, assuming the algebraicity of the K\"unneth components of the
diagonal in the case of Murre's conjecture (see [21], [29], [32]).
Note that the conclusion of Theorem~(0.3) does not imply the
surjectivity of (0.2) for complex algebraic varieties (by taking a
model), because the Leray spectral sequence for a non proper morphism
does not necessarily degenerate at $E_{2}$ (see also [20], [24]).
However, in the case $p=n=2$ (which is related to Voisin's
conjecture), any counter example to the surjectivity does not seem
to be known even over the complex number field.

In this paper we verify the surjectivity in some cases (see
Theorems~(5.5), (5.8) and (5.11)):

\medskip\noindent
{\bf 0.4.~Theorem.} {\it Let $X$ be the complement of a union of
sufficiently general hypersurfaces in a smooth projective variety
over a subfield $k$ of $\C$ \(more precisely, see
Theorems~{\rm (5.5), (5.8)} and {\rm (5.11)).}
In case $k$ is not algebraically closed, replace the target of
{\rm (0.2)} by the category of mixed Hodge structures with
$k$-structure \(defined by using de Rham cohomology\).
Then the \(modified\) morphism {\rm (0.2)} is surjective if $p=n$ or
$p=n+1\le\dim X-2$.
}

\medskip
Note that the converse of Theorem~(0.3) is not true in general
(i.e. the surjectivity of (0.2) does not imply the injectivity of
(0.1), see (5.12)), although it holds in some case where $n=0$
in (0.1) and $n=1$ in (0.2) as treated in [20] (here it is not
necessary to assume $X$ proper for this assertion).
We would need a stronger condition on the surjectivity to show the
injectivity in general.

I would like to thank the referee for useful comments.

\bigskip\bigskip
\centerline{{\bf 1.~Deligne cohomology of mixed sheaves}}

\bigskip\noindent
{\bf 1.1.~Mixed sheaves.}
In this paper $k$ is a subfield of $\C$.
For a $k$-variety $X$, let $\M(X)$ be a category of mixed sheaves
such that the $\M(X)$ satisfy the axioms of mixed sheaves in [30].
More precisely, the $\M(X)$ should be stable by standard functors
like dual, external products, open pull-backs, and the
cohomological direct images by affine morphisms, and they should
satisfy certain compatibility conditions;
then their derived categories $D^{b}\M(X)$ are stable by the
standard functors like direct images and pull-backs, etc.
We also assume that the weight filtration $W$ is defined in $\M(X)$,
and the graded quotients are are polarizable
(e.g. polarizations of Hodge modules are defined over $k$)
so that the pure objects are {\it semisimple},
see loc.~cit.\ for details.

In this paper we assume there exists a natural forgetful functor
$$
\M(X)\to\MHM(X_{\C})
\leqno(1.1.1)
$$
in a compatible way with the above standard functors,
where $\MHM(X_{\C})$ is the category of mixed Hodge modules [28] on
$X_{\C}=X\otimes_{k}\C$.
This condition implies for example that any morphism of $\M(X)$ is
{\it strictly compatible} with $W$.

In case $k$ is not algebraically closed, we assume further a natural
factorization
$$
\M(X)\to\MHM(X)\to\MHM(X_{\C})
\leqno(1.1.2)
$$
in a compatible way with the above standard functors,
where $\MHM(X)$ is the category of mixed Hodge modules on
$X\otimes_{k}\C$ whose underlying filtered $\DD$-module is defined
over $X/k$.

The reader may assume $\M(X)=\MHM(X_{\C})$ if $k$ is
algebraically closed, and $\M(X)=\MHM(X)$ otherwise.
Depending on the purpose, he may also assume some more additional
structure, e.g. systems of realizations $\M_{SR}(X)$ ([30], [32])
which was constructed in [12], [13], [20], etc.\ in the case
$X=\Spec k$.

\medskip\noindent
{\bf 1.2.~Deligne cohomology.} Set $\M(k)=\M(\Spec k)$.
We denote by $\Q$ the constant object in $\M(k)$.
For a $k$-variety $X$, let
$$
\Q_{X}={a}_{X}^{*}\Q,\quad\D_{X}={a}_{X}^{!}
\Q\quad\text{in}\,\,\,\, D^{b}\M(X),
\leqno(1.2.1)
$$
where $a_{X}:X\to\Spec k$ is the structure morphism.
We can define the Tate twist $(m)$ for $m\in\Z$ by using the
cohomology of the projective space $\P^1$ (see e.g. [30]).
We define
$$
H^j(X,\Q)=H^j(a_{X})_*\Q_{X},\quad H_j^{\BM}(X,\Q) =
H^{-j}(a_{X})_*\D_{X}\quad\text{in}\,\,\,\,\M(k),
$$
where $H^j$ is the usual cohomology functor.
They will be denoted by $H^j(X/k,\Q)$ and $H_j^{\BM}(X/k,\Q)$
respectively when we have to specify the ground field $k$ explicitly.

We define analogues of Deligne cohomology and homology by
$$
\aligned H_{\DD}^j(X,\Q(i)) &=\Ext^j(\Q,(a_{X})_*\Q_{X}(i)),
\\ H_j^{\DD}(X,\Q(i)) &=\Ext^{-j}(\Q,(a_{X})_*\D_{X}(-i)),
\endaligned
$$
where $\Ext$ is taken in $D^{b}\M(k)$.

Let $n=\dim X$.
Then we have a canonical morphism
$$
\Q_{X}\to\D_{X}(-n)[-2n]
\leqno(1.2.2)
$$
using the adjunction relation
$$
\Hom(\Q_{X},\D_{X}(-n)[-2n])=\Hom(\Q, H_{2n}^{\BM}(X,\Q)(-n)),
\leqno(1.2.3)
$$
because $H_j^{\BM}(X,\Q)=0$ for $j>2n$ and $H_{2n}^{\BM}(X,\Q)$
is naturally isomorphic to a direct sum of $\Q(n)$ by restricting
to the smooth part of each irreducible component).
If $X$ is smooth and equidimensional, (1.2.2) induces isomorphisms
$$
\aligned
\Q_{X}(n)[n] &=\D_{X}[-n]\quad\text{in}\,\,\,\,\M(X),
\\ H_{\DD}^{2n-j}(X,\Q(n-i)) &=H_j^{\DD}(X,\Q(i)).
\endaligned
\leqno(1.2.4)
$$

Using semisimplicity of pure objects as in [28], 4.5.3, we get
$$
H^{\DD}_j(X,\Q(i))=\Ext^{-j}(\Q,(a_{X})_*\D_{X}(-i))=0\quad
\text{for}\,\,\,j<2i
\leqno(1.2.5)
$$
In the case $X$ is smooth this means
$$
H_{\DD}^j(X,\Q(i))=\Ext^j(\Q,(a_{X})_*\Q_{X}(i))=0\quad
\text{for}\,\,\,j>2i.
$$

\medskip\noindent
{\bf 1.3.~Canonical filtration.}
Let $F_{\tau}$ be a decreasing filtration on $H_{\DD}^j(X,\Q(i))$,
$H_j^{\DD}(X,\Q(i))$ induced by the canonical truncation $\tau$
(see [11]) on $(a_{X})_*\Q_{X}(i)$ and $(a_{X})_*\D_{X}(-i)$.
We shift the filtration so that
$$
\Gr_{F_{\tau}}^aH_{\DD}^j(X,\Q(i)),\quad
\Gr_{F_{\tau}}^aH_j^{\DD}(X,\Q(i))
$$
are respectively subquotients of
$$\Ext^a(\Q, H^{j-a}(X,\Q)(i)),\quad
\Ext^a(\Q, H_{j+a}^{\BM}(X,\Q)(-i)),$$
via the spectral sequences
$$
\aligned
E_2^{p,q}=\Ext^p(\Q, H^q(X,\Q)(i))&\Rightarrow
H_{\DD}^{p+q}(X,\Q(i))
\\
E_2^{p,q}=\Ext^p(\Q, H_{-q}^{\BM}(X,\Q)(-i))&\Rightarrow
H_{-p-q}^{\DD}(X,\Q(i)).
\endaligned
$$
These spectral sequences are associated to the truncation
$\tau_{\le j}$ on $(a_{X})_*\Q_{X}(i)$ and $(a_{X})_*\D_{X}(-i)$
by renumbering the $E_r^{p,q}$ as in [11].

For $a=0,1$, we have canonical injections
$$
\aligned
\Gr_{F_{\tau}}^aH_{\DD}^j(X,\Q(i))&\to\Ext^a(\Q,H^{j-a}(X,\Q)(i)),
\\
\Gr_{F_{\tau}}^aH_j^{\DD}(X,\Q(i))&\to\Ext^a(\Q,H_{j+a}^{\BM}
(X,\Q)(-i)),
\endaligned
\leqno(1.3.1)
$$
since $E_2^{p,q}=0$ for $p<0$ by vanishing of negative extensions.

If $X$ is smooth proper, the filtration $F_{\tau}$ splits by a
variant of the decomposition theorem [30], and (1.3.1) induces
isomorphisms.
The surjectivity of (1.3.1) is not clear except the case
where $\M(k)$ is the category of graded-polarizable mixed
$\Q$-Hodge structures so that higher extension groups $\Ext^i$
vanish for $i>1$ (and hence $E_2^{p,q}=0$ for $p>1$)
as a corollary of Carlson's formula [10].

\medskip\noindent
{\bf 1.4.~Hodge-type conjecture for Chow groups.}
Let $X$ be a smooth proper variety over $k$.
Then the cycle map induces
$$
\CH^p(X)_{\Q}\to\Gr_{F_{\tau}}^{0}H_{\DD}^{2p}(X,\Q(p))=
\Hom_{\M(k)}(\Q, H^{2p}(X,\Q)(p)).
\leqno(1.4.1)
$$
The $\M$-Hodge conjecture means the surjectivity of (1.4.1).
If $k=\C$ and $\M(\C)$ is the category of graded-polarizable
mixed Hodge structures, this is the usual Hodge conjecture.

\medskip\noindent
{\bf 1.5.~Remarks.} (i) Let $U$ be a smooth variety over $k$, and $X$
a smooth compactification.
Then the pull-back
$$
\Hom_{\M(k)}(\Q, H^{2p}(X,\Q)(p))\to
\Hom_{\M(k)}(\Q, H^{2p}(U,\Q)(p))
$$
is surjective by the weight spectral sequence as in [11].
So the $\M$-Hodge conjecture for $U$ can be reduced to that for $X$.

\medskip
(ii) The $\M$-Hodge conjecture for $X$ can be reduced to
the usual Hodge conjecture for $X_{\C}$.
Indeed, a cycle on $X_{\C}$ is defined over a finitely generated
$k$-subalgebra $R$ of $\C$. (This is called a {\it spreading out}.)
We can restrict it to the fiber over a closed point of $\Spec R$.
If $k$ is not algebraically closed, then a cycle $\zeta$ is defined
over a finite Galois extension $k'$ of $k$, and its cycle class in
the de Rham cohomology of $X_{k'}/k'$ is invariant by the action of
$G=\Gal(k'/k)$ using the factorization (1.1.2).
Then we can replace $\zeta$ with
$$
|G|^{-1}\,\msum_{g\in G}\,g^*\zeta,
$$
which is defined over $k$ and whose cycle class is not changed,
see Remark~(iii) below.

\medskip
(iii) For a smooth $k$-variety $X$, let $H^j_{\DR}(X/k)$ denote
the de Rham cohomology of $X/k$.
For a Galois extension $k\subset k'$, set $X_{k'}=X\otimes_kk'$.
Then the Galois group $G:=\Gal(k'/k)$ acts on $X_{k'}\to k'$
and on $H^j_{\DR}(X_{k'}/k')$.
We have a canonical isomorphism
$$
H^j_{\DR}(X_{k'}/k')=H^j_{\DR}(X/k)\otimes_kk',
$$
and the above action is identified with the action associated with
the tensor product with $k'$ over $k$.
Moreover, the cycle map
$$
\CH^j(X_{k'})_{\Q}\to H^{2j}_{\DR}(X_{k'}/k')(j)
$$
is compatible with the action of $G$.

\medskip
(iv) If $X=\Spec K$ with $K$ a finite extension of $k$, then
$$
\Hom_{\MHM(k)}(\Q,H^0(\Spec K/k,\Q))=\Q.
$$
Indeed, we may assume $K=k[t]/(f(t))$ with $f(t)$ irreducible of
degree $d$ over $k$.
Let $\alpha_i$ be the roots of $f(t)$ in $\C$. Then
$$
K\otimes_k\C=\C[t]/(f(t))=\mopl_{i=1}^d\,\C[t]/(t-\alpha_i)=
\mopl_{i=1}^d\,\C.
$$
Moreover, we have in this case
$$
H^0_{\DR}(\Spec K/k)=K.
$$

Consider any morphism in $\MHM(k)$
$$
u:\Q\to H^0(\Spec K/k,\Q)).
$$
Its de Rham part is given by a polynomial $g(t)\in k[t]/(f(t))=K$,
where $\deg g< d$.
Its base change by $k\to\C$ is then given by
$(g(\alpha_i))_{1\le i\le d}$ using the above calculation.
So we get $g(\alpha_i)\in\Q$, since they are the de Rham part of
the morphism
$$
u\otimes_k\C:\Q\to H^0(\Spec K\otimes_k\C/\C,\Q).
$$
This implies that $g(t)\in k$, i.e. $\deg g(t)=0$, since the
$\alpha_i^m\,(0\le m<d)$ are linearly independent over $k$ for each
$i$.

\medskip
(v) Let $X$ be an irreducible projective variety over $k$.
Set $K=k(X)\cap\ok$ in $\overline{k(X)}$.
If $X$ is normal, then $K\in\Gamma(X,\O_X)$, and $X$ is
absolutely irreducible over $K$. So we get in this case
$$
H^0(X/k,\Q)=H^0(\Spec K/k,\Q)\quad\hbox{in}\,\,\,\MHM(k),
$$
and we have by Remark (iv) above
$$
\Hom_{\MHM(k)}(\Q,H^0(X/k,\Q))=\Q.
$$
The last assertion holds without assuming $X$ normal, since we have
the injective morphisms
$$
H^0(\Spec k/k,\Q)\to H^0(X/k,\Q)\to H^0(\tX/k,\Q),
$$
where $\tX$ is the normalization of $X$.

\bigskip\bigskip
\centerline{{\bf 2.~Cycle map of higher Chow groups}}

\bigskip\noindent
{\bf 2.1.~Higher Chow groups} ([8]).
Let $\Delta^n=\Spec (\C[t_0,\dots, t_n]/ (\sum t_i-1))$.
For a subset $I$ of $\{0,\dots, n\}$, let
${\Delta}_I^n=\{t_i=0\,(i\in I)\}\subset\Delta^n$.
We have an inclusion $\iota_i :\Delta^{n-1}\to\Delta^n$ whose
image is ${\Delta}_{\{i\}}^n$ for $0\le i\le n$.

Let $X$ be an equidimensional variety.
We define $\ZZ^p(X,n)$ to be the free abelian group with generators
the irreducible closed subvarieties $Z$ of $X\times\Delta^n$ of
codimension $p$, intersecting al the faces $X\times{\Delta}_I^n$
of $X\times\Delta^n$ properly (i.e. $\dim (Z\cap X\times
{\Delta}_I^n)=\dim Z-|I|)$, see [8].
We have face maps
$$
\partial_i :\ZZ^p(X,n)\to\ZZ^p(X,n-1),
$$
induced by $\iota_i$.
Let $\partial=\sum (-1)^i\partial_i$.
Then $\partial^{2}=0$, and $\CH^p(X,n)$ is defined to be the
homology of the complex (i.e. $\Ker\,\partial /\Im\,\partial)$
which is a subquotient of $\ZZ^p(X,n)$.
Let $\ZZ^p(X,\csbull)'$ be the subcomplex of $\ZZ^p(X,\csbull)$
defined by
$$
\ZZ^p(X,n)'=\mcap_{0\le i<n}\,\Ker(\partial_i :
\ZZ^p(X,n)\to\ZZ^p(X,n-1)).
$$
Then the inclusion induces a quasi-isomorphism
$$
\ZZ^p(X,\csbull)'\to\ZZ^p(X,\csbull)
\leqno(2.1.1)
$$
by [8] (see also [31], 2.1).

\medskip\noindent
{\bf 2.2.~Filtration G.} Set
$$
S_n=\mcup_{0\le i\le n}\,{\Delta}_{\{i\}}^n,\,\,\,
S'_n=\mcup_{0\le i<n}\,{\Delta}_{\{i\}}^n,\,\,\,
U_n=\Delta^n\setminus S_n,\,\,\,
U'_n=\Delta^n\setminus S'_n,
$$
with the inclusions $i_n:S_n\to\Delta^n$, $i'_n:S'_n\to\Delta^n$,
$j_n:U_n\to\Delta^n$, $j'_n:U'_n\to\Delta^n$.
(These morphisms will sometimes denote also the base change of them.)
We have a short exact sequence in $\M(\Delta^n)$
$$
0\to (j_{n-1})_!\Q_{U_{n-1}}[n-1]\to (j_n)_!
\Q_{U_n}[n]\to (j'_n)_!\Q_{U'_n}[n]\to 0,
\leqno(2.2.1)
$$
where $\Delta^{n-1}$ is identified with ${\Delta}_{\{n\}}^n$, and the
direct images by closed embeddings are omitted to simplify the
notation.
This gives an increasing filtration $G$ on $(j_n)_!\Q_{U_n}[n]$
such that
$$
\Gr_{m}^{G}(j_n)_!\Q_{U_n}[n]=(j'_{m})_!
\Q_{U'_{m}}[m]\quad\text{for}\,\,\,0<m\le n,
\leqno(2.2.2)
$$
and $\Gr_0^{G}(j_n)_!\Q_{U_n}[n]=\Q_{\Spec k}$, where $\Spec k$
is identified with the closed point of $\Delta^n$ defined by $t_i=0$
for $i\ne 0$.
Since $U'_n=(\G_{m})^n$, we have
$$
(a_{\Delta^n})_*(j'_n)_!\Q_{U'_n}=0,\quad
(a_{\Delta^n})_*(j_n)_!\Q_{U_n} =\Q_{\Spec k}[-n].
\leqno(2.2.3)
$$

\medskip\noindent
{\bf 2.3.~Cycle map.}
Let $\zeta=\msum_j\,a_j[Z_j]\in\ZZ^p(X,n)'$.
Put $Z=\supp\zeta\,(:=\mcup_j\,Z_j)$, $d=\dim X-p$, and
$d'=\dim Z=d+n$.
We define
$$
u_{\zeta}\in\Ext^{-2d'}(\Q_{Z},\D_{X\times\Delta^n}(-d'))
$$
to be the composition of morphisms
$$
\Q_{Z}\to\mopl_j\,\Q_{Z_j}\to\mopl_j\,
\D_{Z_j}(-d')[-2d']\to\D_{X\times\Delta^n} (-d')[-2d'],
$$
where the second morphism is given by the sum of the canonical
morphisms (1.2.2) multiplied by $a_j$, and the other morphisms are
canonical ones.
Consider the composition of $u_{\zeta}$ with
$$
\D_{X\times\Delta^n}\to (\iota_i)_*{\iota}_i^{*}
\D_{X\times\Delta^n}=\D_{X\times\Delta^{n-1}}(1)[2]
$$
for $i\ne n$, which is induced by
${\iota}_i^{*}:\Q_{\Delta^n}\to\Q_{\Delta^{n-1}}$.
(Here the base change of $\iota_i$ is also denoted by it.)
Let $Z_{(i)}=Z\cap X\times {\Delta}_{\{i\}}^n$.
Then the composition vanishes by the assumption that
$\zeta\in\Ker\,\partial_i$.
Indeed, it is identified with an element of
$$
\aligned
\Ext^{2-2d'}(\Q_{Z_{(i)}},\D_{Z_{(i)}}(1-d'))
&=\Hom(\Q,H_{2d'-2}^{\BM}(Z_{(i)},\Q)(1-d'))
\\
&=\mopl_j\,\Q
\endaligned
$$
by using the adjunction relations for the inclusion
$\iota_i:X\times\Delta^{n-1}\to X\times
\Delta^n,$ and is given by the intersection multiplicity of $\zeta$
with $X\times {\Delta}_{\{i\}}^n$ at each irreducible component of
$Z_{(i)}$, because the cycle map is compatible with the pull-back
for a closed immersion of a locally principal divisor
(see e.g. [29], II).
Combined with vanishing of negative extensions, this implies
vanishing of the composition of $u_{\zeta}$ with
$$
\D_{X\times\Delta^n}\to (i'_n)_*(i'_n)^{*}
\D_{X\times\Delta^n},
$$
because $(i'_n)^{*}\D_{\Delta^n}=\Q_{S'_n}(n)[2n]$ and the
graded-pieces of the weight filtration on
$\Q_{S'_n}[n-1]\in\M(S'_n)$ are constant sheaves supported on
intersections of faces.
So $u_{\zeta}$ is uniquely lifted to
$$
u'_{\zeta}\in\Ext^{-2d'}(\Q_{Z}, (j'_n)_!\D_{X\times U'_n}(-d'))
$$
by using the long exact sequence together with vanishing of
negative extensions.

If furthermore $\zeta\in\Ker\,\partial_n$, we see that $u'_{\zeta}$
is uniquely lifted to
$$
u''_{\zeta}\in\Ext^{-2d'}(\Q_{Z}, (j_n)_!\D_{X\times U_n}(-d'))
$$
by a similar argument.
Taking the composition with $\Q_{X\times\Delta^n}\to\Q_{Z}$ we get
$$
\aligned v_{\zeta} &\in\Ext^{-2d'}(\Q_{X\times\Delta^n},
(j_n)_!\D_{X\times U_n}(-d'))
\\
&=\Ext^{-2d'}(\Q,(a_{X\times\Delta^n})_* (j_n)_!\D_{X\times U_n}
(-d')),
\endaligned
$$
and it defines $cl(\zeta)\in H_{2d+n}^{\DD}(X,\Q(d))$, because
$(a_{\Delta^n})_*(j_n)_!\D_{U_n}=\Q(n)[n]$.

This construction defines the cycle map
$$
cl :\CH^p(X,n)\to H_{2d+n}^{\DD}(X,\Q(d)).
\leqno(2.3.1)
$$
Indeed, if $\zeta$ belongs to the image of the differential of
$\ZZ^p(X,\csbull)'$, then $v_{\zeta}$ comes from
$$
\Ext^{-2d'-2}(\Q, (a_{X\times\Delta^{n+1}})_* (j_{n+1})_!
\D_{X\times U_{n+1}}(-d'-1)),
$$
(which vanishes by (2.2.3)) by using the compatibility of the cycle
map and the pull-back for $X\times\Delta^n\to X\times\Delta^{n+1}$.

If $X$ is smooth, (2.3.1) gives
$$
cl :\CH^p(X,n)\to H_{\DD}^{2p-n}(X,\Q(p)).
\leqno(2.3.2)
$$
If $n=0$ and $X$ is smooth proper, this induces the Abel-Jacobi map
[18] using Carlson's formula [10].
We can show that the cycle map is compatible with the pushforward by
a proper morphism and the pull-back by a morphism of smooth
quasi-projective varieties, see [31].
So it is also compatible with the action of algebraic correspondences.

Let $F_{\tau}$ denote also the induced filtration on $\CH^p(X,n)$
by the cycle map, see (1.3).
We have the induced morphisms
$$
\Gr_{F_{\tau}}^icl :\Gr_{F_{\tau}}^i\CH^p(X,n)
\to\Gr_{F_{\tau}}^iH_{\DD}^{2p-n}(X,\Q(p)).
\leqno(2.3.3)
$$
Note that the target is $\Ext^i(\Q, H^{2p-n-i}(X,\Q)(p))$ if $X$
is smooth proper, see (1.3).

\medskip\noindent
{\bf 2.4.~Lemma.}
{\it Let $X, Y$ be smooth projective varieties over $k$.
If the $\M$-Hodge conjecture holds for cycles of codimension
$\dim Y-j$ in $X\times_{k} Y$, then for
$\zeta\in\CH^{j+\dim X}(X\times_{k} Y)_{\Q}$, we have
$$
\zeta_*(\Im\,\Gr_{F_{\tau}}^icl) =
\Im\,\Gr_{F_{\tau}}^icl\cap\Im\,\zeta_*\quad
\text{in}\,\,\,\,\Ext^i(\Q,H^{2p-n-i}(Y,\Q)(p)),
\leqno(2.4.1)
$$
where
$$
\zeta_*:\Ext^i(\Q,H^{q-2j}(X,\Q)(p-j))\to\Ext^i(\Q,H^{q}(Y,\Q)(p))
$$
is induced by $\zeta$, and $\Gr_{F_{\tau}}^icl$ is as in
{\rm (2.3.3).} }

\medskip\noindent
{\it Proof.} By assumption, there exists
$\zeta'\in\CH^{\dim Y-j}(Y\times_{k} X)_{\Q}$ such that
$$
\zeta'_*:H^{q}(Y,\Q)\to H^{q-2j}(X,\Q)(-j)
$$
vanishes for $q\ne 2p-n-i$, and the restriction of
$\zeta_*\scirc\zeta'_*$ to $\Im\,\zeta_*\subset H^{q}(Y,\Q)$ is the
identity for $q=2p-n-i$, because $H^{q-2j}(X,\Q)$ and $H^{q}(Y,\Q)$
are semi\-simple.
So the assertion follows from the compatibility of (2.3.3) with the
action of correspondences.

\bigskip\bigskip
\centerline{{\bf 3.~Case of varieties with normal crossings}}

\bigskip\noindent
{\bf 3.1.~Variant for higher Chow groups.} We say that $Y$ is a
variety with normal crossings if $Y$ is equidimensional, and is
locally isomorphic to a divisor with normal crossings.
In this paper we assume that the irreducible components
$Y_i\,(1\le i\le r)$ are smooth.
For a subset $I$ of $\{1,\dots, r\}$, we define
$Y_I=\mcap_{i\in I}\,Y_i$.

Consider the double complex
$$
\to\bigoplus_{|I|=i+1}\ZZ^{p-i}(Y_I,\csbull)
\to\bigoplus_{|I|=i}\ZZ^{p-i+1}(Y_I,\csbull)
\to\,\,\,\cdots\,\,\,
\to\bigoplus_{|I|=1}\ZZ^p(Y_I,\csbull)\to 0,
$$
where the differential is given by the alternating sum as usual,
and $\mopl_{|I|=1}\,\ZZ^p(Y_I,\csbull)$ is placed at the degree
zero.
Let $\tcZ^p(Y,\csbull)$ be its total complex, and similarly for
$\tcZ^p(Y,\csbull)'$ (with $\ZZ^p(Y_I,\csbull)$ replaced by
$\ZZ^p(Y_I,\csbull)')$.
We have a canonical quasi-isomorphism
$\tcZ^p(Y,\csbull)'\to\tcZ^p(Y,\csbull)$ by (2.1.1).

We can verify that the canonical morphisms
$$
\tcZ^p(Y,\csbull)\to\ZZ^p(Y,\csbull),\quad
\tcZ^p(Y,\csbull)'\to\ZZ^p(Y,\csbull)'
\leqno(3.1.1)
$$
are quasi-isomorphisms.
(Indeed, for an irreducible subvariety $Z$ of $Y\times\Delta^n$
intersecting all the faces $Y\times {\Delta}_I^n$ properly,
let $J$ be the subset of $\{1,\dots, r\}$ such that
$Y_i\times\Delta^n$ contains $Z$ for $i\in J$.
Then $Z$ defines an element of $\ZZ^{p-|I|+1}(Y_I,n)$ for
$I\subset J$.
This gives an acyclic complex if we define an augmented complex
by adding the term for $Y_{\emptyset}=Y$.
So we get the first quasi-isomorphism.
Then the second follows from (2.1.1).)

Let $W$ be an increasing filtration on $\tcZ^p(Y,\csbull)'$ such
that $\Gr_i^{W}\tcZ^p(Y,\csbull)'$ consists of
$\mopl_{|I|=i+1}\,\ZZ^{p-i}(Y_I,\csbull)$ for $i\ge 0$ and
$W_{-1}\tcZ^p(Y,\csbull)=0$.
We define
$$
\tcZ^p(Y^{[a,b]},\csbull)'=
W_{b}\tcZ^p(Y,\csbull)'/W_{a-1}\tcZ^p(Y,\csbull)'.
$$
Let $\CH^p(Y^{[a,b]},n)$ be its homology.
We denote them by $\tcZ^p(Y^{[a]},\csbull)',\CH^p(Y^{[a]},n)$
if $a=b$, and by $\tcZ^p(Y^{\ge a},\csbull)',\CH^p(Y^{\ge a},n)$
if $b=+\infty$.
Then we have a canonical long exact sequence
$$
\CH^p(Y^{[a]},n)
\to\CH^p(Y^{\ge a},n)
\to\CH^p(Y^{\ge a+1},n)
\to\CH^p(Y^{[a]},n-1)
\leqno(3.1.2)
$$
Note that we have by definition
$$
\CH^p(Y^{[a]},n)=\mopl_{|I|=a+1}\,\CH^{p-a}(Y_I, n-a).
\leqno(3.1.3)
$$

\medskip\noindent
{\bf 3.2.~Variant for Deligne cohomology.} Let $W$ be the weight
filtration on $\D_{Y}\in\M(Y)[\dim Y]$ (which is a mixed sheaf
shifted by $\dim Y)$ so that
$$
\Gr_{a}^{W}\D_{Y}=\mopl_{|I|=a+1}\,\D_{Y_I}[a]
\in\M(Y)[\dim Y].
$$
Let $\D_{Y}^{[a,b]}=W_{b}\D_{Y}/W_{a-1}\D_{Y}$, and
$$
H_j^{\DD}(Y^{[a,b]},\Q(i))=
\Ext^{-j}(\Q,(a_{Y})_*\D_{Y}^{[a,b]}(-i)).
$$
We denote it by $H_j^{\DD}(Y^{[a]},\Q(i))$ if $a=b$, and by
$H_j^{\DD}(Y^{\ge a},\Q(i))$ if $b=+\infty$.
Then we have a long exact sequence
$$
H_j^{\DD}(Y^{[a]},\Q(i))
\to H_j^{\DD}(Y^{\ge a},\Q(i))
\to H_j^{\DD}(Y^{\ge a+1},\Q(i))
\to H_{j-1}^{\DD}(Y^{[a]},\Q(i)).
\leqno(3.2.1)
$$
Note that we have by definition
$$
H_j^{\DD}(Y^{[a]},\Q(i))=
\mopl_{|I|=a+1}\,H_{j-a}^{\DD}(Y_I,\Q(i)).
\leqno(3.2.2)
$$

The filtration $W$ induces naturally a filtration $W$ on
$\D_{Y\times U_n}$ which is the external product of $\D_{Y}$ and
$\D_{U_n}$.
Let
$$
\D_{Y\times {U}_n}^{\ge a}=\D_{Y\times U_n}/W_{a-1}\D_{Y\times U_n}.
$$
We define an increasing filtration $G$ on
$$
(j_n)_!\D_{Y\times {U}_n}^{\ge a}\in\M(Y\times
\Delta^n)[n +\dim Y]
$$
to be the convolution (see [5]) of $W$ on $\D_{Y}^{\ge a}$ and
$G$ on $(j_n)_!\D_{U_n}$, where the latter is induced by $G$ on
$$
(j_n)_!\Q_{U_n}[n]=(j_n)_!\D_{U_n} (-n)[-n].
$$
(Note that $(j_n)_!\D_{Y\times {U}_n}^{\ge a}$ is the external
product of $\D_{Y}^{\ge a}$ and $(j_n)_!\D_{U_n}$, because the
direct image is compatible with external product.)
Then $\Gr_{m}^{G}(j_n)_!\D_{Y\times {U}_n}^{\ge a}$ is the direct
sum of the external products of
$$
\mopl_{|I|=i+1}\,\D_{Y_I}[i]\quad\text{and}\quad
(j'_{m-i})_!\D_{U'_{m-i}}(n-m+i)[n-m+i]\quad\text{for}\,\,\,i\ge a.
$$

\medskip\noindent
{\bf 3.3.~Proposition.} {\it We have cycle maps
$$
cl :\CH^p(Y^{\ge a},n)\to H_{2d+n}^{\DD}(Y^{\ge a},
\Q(d))
\leqno(3.3.1)
$$
which induce a morphism of the long exact sequence {\rm (3.1.2)} to
{\rm (3.2.1)}, where $d=\dim Y-p$, and we use {\rm (2.3.1)} for
$\CH^p(Y^{[a]},n)$.
Furthermore {\rm (3.3.1)} is identified with {\rm (2.3.1)} if $a=0$.
}

\medskip\noindent
{\it Proof.}
Let $\zeta=\msum_{i\ge a}\,\zeta_i\in {\tcZ}^p(Y^{\ge a},n)'$ with
$\zeta_i\in\mopl_{|I|=i+1}\,\ZZ^{p-i}(Y_I,n) $.
Put $Z=\supp\zeta, Z_i=\supp\zeta_i$ so that $Z=\mcup_{i\ge a}\,Z_i$.
Let $d, d'$ be as in (2.3).
Then $\zeta$ defines
$$
u'_{\zeta}\in\Ext^{-2d'}(\Q_{Z},
\Gr_n^{G}(j_n)_!\D_{Y\times {U}_n}^{\ge a}(-d'))
$$
by an argument similar to (2.3).
If $\zeta$ is annihilated by the differential of
${\tcZ}^p(Y^{\ge a},n)'$, then it is uniquely lifted to
$$
u'_{\zeta}\in\Ext^{-2d'}(\Q_{Z}, G_n(j_n)_!
\D_{Y\times {U}_n}^{\ge a}(-d')),
$$
by using the compatibility of the cycle map with the pull-back by a
closed immersion of a principal divisor.
Taking the composition with natural morphisms, we get
$$
\aligned
v_{\zeta} &\in\Ext^{-2d'}(\Q_{Y\times\Delta^n},
(j_n)_!\D_{Y\times {U}_n}^{\ge a}(-d'))
\\ &=\Ext^{-2d'}(\Q, (a_{Y\times\Delta^n})_*(j_n)_!
\D_{Y\times {U}_n}^{\ge a}(-d')),
\endaligned
$$
and it gives $cl(\zeta)\in H_{2d+n}^{\DD}(Y^{\ge a},\Q(d))$.

This defines a well-defined cycle map
$$
cl :\CH^p(Y^{\ge a},n)\to H_{2d+n}^{\DD}(Y^{\ge a},\Q(d)).
$$
Indeed, if $\zeta$ belongs to the coboundary, then the image of
$u'_{\zeta}$ in
$$
\Ext^{-2d'}(\Q_{Y\times\Delta^{n+1}},G_{n+1}(j_{n+1})_!
\D_{Y\times {U}_{n+1}}^{\ge a}(-d'))
$$
vanishes by the long exact sequence associated with
$$
0\to G_n\to G_{n+1}\to\Gr_{n+1}^{G}\to 0.
$$
The remaining assertions follow from the construction easily.
This finishes the proof of Proposition~(3.3).

\bigskip\bigskip
\centerline{{\bf 4.~Proof of Theorem~(0.3)}}

\bigskip\noindent
In this section we prove Theorem~(0.3) by showing Theorems~(4.1) and
(4.2) below.
Here the condition $k=\oQ$ in Theorem~(0.3) does not appear
explicitly.
This is implicitly used in the assumption on the injectivity of
(4.1.1).

\medskip\noindent
{\bf 4.1.~Theorem.}
{\it Let $U$ be a smooth quasi-projective variety over $k$, and $X$
a smooth projective compactification of $X$ such that the complement
$Y$ is a divisor with normal crossings whose irreducible components
$Y_i$ are smooth.
Let $Y_I$ be as in {\rm (3.1)}, and put $Y_{\emptyset}=X$ for
$I=\emptyset$.
Let $p,n$ be integers, and set $d=\dim X-p$.
Assume that the cycle maps
$$
\CH^{p-|I|}(Y_I,n-|I|)_{\Q}\to
H_{\DD}^{2p-n-|I|} (Y_I,\Q(p-|I|))/F_{\tau}^{2}
\leqno(4.1.1)
$$
are injective and the $\M$-Hodge conjecture for $Y_I$ and
$Y_I\times Y_{I\cup\{i\}}$ is true for any $I$ and $i\notin I$
\(including the case $I=\emptyset)$.
Then the following cycle maps are surjective for any $a\ge 0$ in the
notation of {\rm (3.1):}}
$$
\CH^p(U,n+1)_{\Q}
\to\Gr_{F_{\tau}}^{0}H_{\DD}^{2p-n-1}(U,\Q(p)),
\leqno(4.1.2)
$$$$
\CH^{p-1}(Y^{\ge a},n)_{\Q}
\to\Gr_{F_{\tau}}^{0}H_{2d+n}^{\DD}(Y^{\ge a},\Q(d)).
\leqno(4.1.3)
$$

\medskip\noindent
{\it Proof.} We first show that the surjectivity of (4.1.2)
assuming that for (4.1.3).
We may assume $n\ge 0$ because the case $n+1=0$ follows from
the $\M$-Hodge conjecture for $X$, see Remark~(1.5)(i).
Note that the target of (4.1.2) vanishes for $n+1<0$ since
$H^{2p-n-1}(U,\Q(p))$ has weights $\ge -n-1$, see also (1.2.5).
Similarly the target of (4.1.3) vanishes for $n<a$ using (3.2.2)
together with a spectral sequence.

Set $q=2d+n$.
By an argument similar to [31], we have a commutative diagram
$$
\CD @>>>\CH^p(U,n+1)_{\Q} @>>>\CH^{p-1}(Y,n)_{\Q} @>>>
\CH^p(X,n)_{\Q} @>>>
\\ @. @VVV @VVV @VVV
\\ @>>> H_{q+1}^{\DD}(U,\Q(d)) @>>> H_q^{\DD}(Y,\Q(d)) @>>>
H_q^{\DD}(X,\Q(d)) @>>>
\endCD
$$
Let $\xi\in H_{q+1}^{\DD}(U,\Q(d))$, and $\xi'$ be its image in
$H_q^{\DD}(Y,\Q(d))$.
By the surjectivity of (4.1.3) for $a=0$, there exists
$\zeta'\in\CH^{p-1}(Y,n)_{\Q}$ such that
$cl(\zeta')-\xi'\in F_{\tau}^1H_q^{\DD}(Y,\Q(d))$.
Let $\zeta''$ be the image of $\zeta'$ in $\CH^p(X,n)_{\Q}$.
Then $cl(\zeta'')$ coincides with the image of $cl(\zeta')-\xi'$
and belongs to $F_{\tau}^1H_q^{\DD}(X,\Q(d))$.

By the weight spectral sequence [11], we see that the morphisms
$Y_i\to Y\to X$ induce an isomorphism
$$
\Im(H_{q+1}^{\BM}(Y,\Q)\to H_{q+1}^{\BM}(X,\Q)) =
\Im(\mopl_i\,H_{q+1}^{\BM}(Y_i,\Q)\to H_{q+1}^{\BM}(X,\Q)),
\leqno(4.1.4)
$$
and the projection
$$
H_{q+1}^{\BM}(Y,\Q)\to\Im(H_{q+1}^{\BM}(Y,\Q)\to H_{q+1}^{\BM}(X,\Q))
$$
splits by semi\-simplicity of $H_{q+1}^{\BM}(Y_i,\Q)\in\M(k)$.
So there exists $\zeta_{1}\in\CH^{p-1}(Y,n)_{\Q}$ by Lemma~(2.4)
(applied to $X$ and $Y_i)$ such that
$cl(\zeta_{1})\in F_{\tau}^1H_q^{\DD}(Y,\Q(d))$ and the images
of $cl(\zeta')$ and $cl(\zeta_{1})$ in
$\Gr_{F_{\tau}}^1H_q^{\DD}(X,\Q(d))$ coincide.
Thus, replacing $\zeta'$ with $\zeta'-\zeta_{1}$, we may assume
$$
cl(\zeta'')\in F_{\tau}^{2}H_q^{\DD}(X,\Q(d)).
$$
But this implies $\zeta''=0$ by the injectivity of (4.1.1), and
$\zeta'$ comes from $\zeta\in\CH^p(U,n+1)_{\Q}$.
So the surjectivity of (4.1.2) is reduced to the injectivity of
$$
\Hom_{\M(k)}(\Q,H_{q+1}^{\BM}(U,\Q))\to
\Hom_{\M(k)}(\Q,H_q^{\BM}(Y,\Q)).
$$
Using the long exact sequence of Borel-Moore homology
$$
\to H_{q+1}^{\BM}(X,\Q)
\to H_{q+1}^{\BM}(U,\Q)
\to H_q^{\BM}(Y,\Q)\to,
$$
together with semi\-simplicity of $H_{q+1}^{\BM}(X,\Q)$, this
injectivity follows from vanishing of
$\Hom_{\M(k)}(\Q,H_{q+1}^{\BM}(X,\Q)(d))$ (where the target is
pure of weight $-n-1<0$).

The proof of the surjectivity of (4.1.3) is by decreasing induction
on $a$, and is similar to the above argument.
Here we use Proposition~(3.3) instead of the above commutative
diagram, and the isomorphism corresponding to (4.1.4) follows from
the morphisms
$$
\D_{Y}^{[a+1]}\to\D_{Y}^{\ge a+1}\to\D_{Y}^{[a]}[1].
$$
Note that the surjectivity of (4.1.3) at the initial stage of the
induction (where $Y^{\ge a}=Y^{[a]}$) follows from the Hodge
conjecture, because the target of (4.1.3) vanishes if $n-a\ne 0$.
This completes the proof of Theorem~(4.1).

\medskip\noindent
{\bf 4.2.~Theorem.} {\it With the notation and the assumptions of
Theorem~{\rm (4.1)}, the following cycle maps are injective for any
$a\ge 0:$}
$$
\CH^p(U,n)_{\Q}\to H_{\DD}^{2p-n}(U,\Q(p))/F_{\tau}^{2},
\leqno(4.2.1)
$$
$$
\CH^{p-1}(Y^{\ge a},n-1)_{\Q}\to
H_{2d+n-1}^{\DD}(Y^{\ge a},\Q(d))/F_{\tau}^{2}.
\leqno(4.2.2)
$$

\medskip\noindent
{\it Proof.} We first show the assertion for (4.2.1) assuming
that for (4.2.2).
Set $q=2d+n$.
Consider the commutative diagram
$$
\CD
\CH^{p-1}(Y,n)_{\Q} @>>>\CH^p(X,n)_{\Q} @>>>
\CH^p(U,n)_{\Q} @>>>\CH^{p-1}(Y,n-1)_{\Q}
\\ @VVV @VVV @VVV @VVV
\\ H_q^{\DD}(Y,\Q(d)) @>>> H_q^{\DD}(X,\Q(d)) @>>>
H_q^{\DD}(U,\Q(d)) @>>> H_{q-1}^{\DD}(Y,\Q(d))
\endCD
$$
Let $\zeta\in\CH^p(U,n)_{\Q}$, and assume
$cl(\zeta)\in F_{\tau}^{2}H_q^{\DD}(U,\Q(d))$.
Then the image of $\zeta$ in $\CH^{p-1}(Y,n-1)_{\Q}$ vanishes by
(4.2.2), and $\zeta$ comes from $\zeta'\in\CH^p(X,n)_{\Q}$.
Let $\xi'=cl(\zeta')$.
Note that
$$
\Gr_{F_{\tau}}^{0}\xi'\in\Hom(\Q,H_q^{\BM}(X,\Q)(d))=0,
$$
unless $n=0$.
Since the image of $\xi'$ in $\Gr_{F_{\tau}}^{0}H_q^{\DD}(U,\Q(d))$
vanishes, it follows from (4.1.4) that $\Gr_{F_{\tau}}^{0}\xi'$
comes from
$$
\xi''\in\mopl_i\,\Gr_{F_{\tau}}^{0}H_q^{\DD} (Y_i,\Q(d)).
$$
Replacing $\zeta'$ if necessary, we may assume
$\Gr_{F_{\tau}}^{0}\xi'=0$ by the Hodge conjecture for $Y_i$.
This implies that $\xi'$ belongs to
$\Im\,cl\cap F_{\tau}^1H_q^{\DD}(X,\Q(d))$, and hence
$\Gr_{F_{\tau}}^1\xi'$ comes from $\CH^{p-1}(Y,n)_{\Q}$ by an
argument similar to the proof of Theorem~(4.2) (using Lemma~(2.4)).
So we may assume $\xi'\in F_{\tau}^{2}H_q^{\DD}(X,\Q(d))$
(replacing $\zeta'$ if necessary), and the assertion for (4.2.1)
follows from the hypothesis on the injectivity of (4.1.1).

The proof of the injectivity of (4.2.2) is similar by decreasing
induction on $a$ using a morphism of (3.1.2) to (3.2.1) (see
Proposition~(3.3)) instead of the above commutative diagram.
This finishes the proof of Theorem~(4.2).

\medskip\noindent
{\bf 4.3.~Proof of Theorem~(0.3).} The first assertion follows from
Theorems~(4.1) and (4.2) together with (1.3).
Then the last assertion is clear by definition of the refined cycle
map using the spreading out of algebraic cycles ([6], [17], [34]).

\medskip\noindent
{\bf 4.4.~Remark.} Let $X$ be a smooth projective variety over $\C$.
If (0.2) is surjective for any open subvarieties of $X$ with $p=n=2$,
then it would imply the injectivity of the higher Abel-Jacobi map of
the indecomposable higher Chow group to the reduced Deligne cohomology
$$
\CH_{\ind}^{2}(X,1)_{\Q}\to
\Ext_{\MHS}^1(\Q,H^{2}(X,\Q)(2))/
\Hdg^1(X)_{\Q}\otimes_{\Q}{\C}^{*},
\leqno(4.4.1)
$$
where $\Hdg^1(X)_{\Q}=\Hom_{\MHS}(\Q,H^{2}(X,\Q)(1))$,
see [3], [23], [26], [27], [31].
The image of (4.4.1) is countable by a rigidity argument
(see [2], [23]), and its injectivity would imply Voisin's
conjecture [35] on the countability of the indecomposable higher
Chow group.
Note that this conjecture can also be reduced to the injectivity
of the refined cycle map (see [32]), and to the hypotheses of
Theorem~(0.3).

\bigskip\bigskip
\centerline{{\bf 5.~Complement of general hypersurfaces}}

\bigskip\noindent
{\bf 5.1.~Moderate singularities.} Let $Z$ be a complex algebraic
variety.
Assume $Z$ is purely $d$-dimensional and $\Q_{Z}[d]$ is a perverse
sheaf (e.g. $Z$ is a divisor on a smooth variety).
Let $W$ be the weight filtration on $\Q_{Z}[d]$, see [5], [28].
We say that $Z$ has only {\it moderate singularities} in this paper,
if
$$
\dim\supp\Gr_i^{W}(\Q_{Z}[d])<i\quad\text{for}\,\,\,i<d.
\leqno(5.1.1)
$$
If $Z$ is a $k$-variety, we say that $Z$ has only moderate
singularities if so is $Z_{\C} :=Z\otimes_{k}\C$.
Condition (5.1.1) is trivially satisfied if $Z$ is smooth, or more
generally, if $Z$ is a $\Q$-homology manifold.
By duality, (5.1.1) is equivalent to
$$
\dim\supp\Gr_{i-d}^{W}(\D_{Z}[-d])<d-i\quad\text{for}\,\,\, i>0,
\leqno(5.1.2)
$$
and implies
$$
\Gr_{2i-2d}^{W}{H}_{2d-i}^{BM}(Z,\Q)=0\quad\text{for}\,\,\,i>0.
\leqno(5.1.3)
$$

If $Z$ is a divisor on a smooth variety $X$ with a local (reduced)
defining equation $f$, then (5.1.1) is equivalent to
$$
\dim\supp\Gr_{d+i}^{W}\varphi_{f,1}\Q_{X}[d])<d-i\quad
\text{for}\,\,\,i>0.
\leqno(5.1.4)
$$
using the short exact sequence
$$
0\to\Q_{Z}[d]\to\psi_{f,1}\Q_{X}[d]\to\varphi_{f,1}\Q_{X}[d]\to 0
$$
together with
$$
\Q_{Z}[d]=\Ker\,N\subset\psi_{f,1}\Q_{X}[d],
$$
because the weight filtration $W$ on $\psi_{f,1}\Q_{X}[d]$ is the
monodromy filtration shifted by $d$.
Here $\psi_{f,1}$ and $\varphi_{f,1}$ denote the unipotent
monodromy part of Deligne's nearby and vanishing cycle functors
$\psi_{f}$ and $\varphi_{f}$ respectively, and $N=\log T_{u}$ with
$T=T_{u}T_{s}$ the Jordan decomposition of the monodromy $T$.
If furthermore $Z$ has only isolated singularities, then (5.1.1)
is equivalent to that the Jordan blocks of the Milnor monodromy
$T$ on the vanishing cohomology for the eigenvalue $1$ have size
$< d$.
In the case $d=1$, the condition is equivalent to the analytic
local irreducibility.

If $Z$ is a divisor on a smooth proper variety $X$, and $U$ is
its complement, then condition (5.1.1) implies
$$
\Gr_{2i+2}^{W}H^{i+1}(U,\Q)=0\quad\text{for}\,\,\,i>0.
\leqno(5.1.5)
$$
Note that (5.1.5) holds also for $i=0$ if $Z$ is irreducible.

If $Z$ is a quasi-projective variety, then condition (5.1.1) is
stable by a generic hypersurface section.
Note that the irreducibility of $Z$ is also stable by a generic
hypersurface section of positive dimension.
(This follows for example from a generalization of the weak
Lefschetz theorem [5] applied to a smooth affine open subvariety.)

\medskip\noindent
{\bf 5.2.~Lemma.}
{\it Let $X$ be a connected smooth complex projective variety,
and $D$ a divisor on $X$ such that $X' :=X\setminus D$ is affine.
Let $Y$ be a smooth hypersurface section of $X$ which intersects
transversely each stratum of a Whitney stratification of $X$
compatible with $D$.
Put $Y'=Y\cap X'$ with the inclusion $i':Y'\to X'$.
Assume that the cycle class of $Y'$ vanishes in $H^{2}(X',\Q)(1)$.
Then the localization sequence induces the exact sequences
$$
0\to H^j(X',\Q)\to H^j(X\setminus Y',\Q)\to H^{j-1}(Y',\Q)(-1)
\to 0.
\leqno(5.2.1)
$$
}

\medskip\noindent
{\it Proof.} It is enough to show vanishing of the Gysin morphism
$$
i'_*:H^j(Y',\Q)\to H^{j+2}(X',\Q)(1)
\leqno(5.2.2)
$$
for any $j$.
Let $d=\dim X$.
Then $H^j(X',\Q)=0$ for $j>d$, and $H^j(Y',\Q)=0$ for $j>d-1$.
The restriction morphism
$$
i^{\prime *}:H^j(X',\Q)\to H^j(Y',\Q)
\leqno(5.2.3)
$$
is an isomorphism for $j\le d-2$ by a generalization of the weak
Lefschetz theorem for perverse sheaves [5] (applied to the perverse
sheaf $\mathbf{R}j'_*\Q_{X'}[d]$ where $j':X'\to X$ denotes the
inclusion).
So it is enough to show vanishing of the composition of (5.2.3)
and (5.2.2).
But this composition coincides with the cup product with the
cohomology class of $Y'$, and it vanishes by hypothesis.
So the assertion follows.

\medskip\noindent
{\bf 5.3.~Hypersurface sections.} Let $X$ be a geometrically
irreducible smooth projective $k$-variety, where $k$ is a
subfield of $\C$.
For a line bundle $L$, we define
$$
\PP_{L}=\Proj(\Sym_{k}\Gamma (X,L)^{\vee}).
$$
(It is the projective space associated with the symmetric algebra
of the dual vector space of $\Gamma (X,L)$ over $k$.) Note that a
$k$-valued point $z$ of $\PP_{L}$ corresponds to a divisor $D$ on
$X$ such that $\O_{X}(D)\simeq L$.

For $0\le i\le m$, let $L_i$ be line bundles, and
$z_i\in\PP_{L_i}(k)$.
Then $z_i$ corresponds to a divisor $D_i$ (which is also denoted
by $D_{z_i}$) on $X$ as above.
Consider the canonical morphism
$$
\mopl_{0\le i\le m}\Q [D_i]\to\CH^1(X)_{\Q},
\leqno(5.3.1)
$$
where the source is a $\Q$-vector space with basis $[D_i]$.
Let $u_{1},\dots, u_{r}$ be a basis of the kernel of (5.3.1) such
that $u_j\in\mopl_{0\le i\le m}\Z [D_i]$.
Then there is a rational function $f_j$ on $X$ such that
$\div f_j=u_j$ (replacing $u_j$ if necessary).
The rational function $f_j$ is identified with an element of
$\CH^1(X',1)$ for an open subvariety $X'$ of $X$ such that $f_j$
has no zeros nor poles on $X'$.
We will denote by $df_j/f_j$ the image of $f_j$ by the cycle map
$$
\CH^1(X',1)\to\Hom_{\M(k)}(\Q, H^1(X',\Q)(1)).
\leqno(5.3.2)
$$
Indeed, the image is expressed by $df_j/f_j$ at the level of de
Rham cohomology.
Note that in the case $k=\C$ and $\M(k)=\MHS$, the morphism
(5.3.2) is surjective with kernel $\C^{*}$, and an element in the
target of (5.3.2) is called an integral logarithmic $1$-form (i.e.
of the form $df/f$ for a rational function $f$) if it comes from
integral cohomology.

We define for $I=\{i_{1},\dots, i_j\}$
$$
df_I/f_I=df_{i_{1}}/f_{i_{1}}\wedge\dots\wedge df_{i_j}/f_{i_j}\in
\Hom_{\M(k)}(\Q, H^j(X',\Q)(j)).
$$

\medskip\noindent
{\bf 5.4.~Generic condition.} With the notation of (5.3) we assume
$L_i$ are very ample for $i>0$.
Let $\PP'$ be the open subvariety of
$\PP :=\prod_{1\le i\le m}\PP_{L_i}$ such that the divisors $D_{z_i}$
corresponding to $z=(z_{1},\dots, z_{m})\in\PP'(\ok)$ are smooth and
intersect each other and also $D_0$ transversely (more precisely,
there exists a Whitney stratification of $X$ compatible with $D_0$
such that the restrictions of $D_{z_i}$ to each stratum form a
divisor with normal crossings).
We consider further the subset $P'_0$ of $\PP'(k)$ consisting of $z$
which satisfies the following condition
(which is closely related to Remark~(1.5)(iv)):

\medskip\noindent
(5.4.1)\,\,\, For any subset $I$ of $\{0,\dots, m\}$ such that
$|I|=n$, the intersection $D_{I,z} :=\mcap_{i\in I}\,D_{z_i}$
consists of one point (i.e. the Galois group $\Gal(\ok/k)$ acts
transitively on the $\ok$-valued points in the intersection).

\medskip
Let $D_I$ be the closed subvariety of $X\times_{k}\PP'$
with the projection $\pi:D_I\to\PP'$ such that the fiber over
$z\in\PP'(\ok)$ is $D_{I,z}$.
Then $D_I$ is irreducible by a monodromy argument, and Hilbert's
irreducibility theorem asserts that $P'_0$ is quite large in the
case $k$ is finitely generated over $\Q$, see [22], [33].

Let $k_0$ be a subfield of $k$ such that $X, L, D_0,\PP$ and $\PP'$
are defined over $k_0$.
We say that $(D_{1},\dots, D_{m})$ is $k_0$-{\it generic}, if the
corresponding $z=(z_{1},\dots, z_{m})\in\PP(k)$ is a generic point
of $\PP$ relative to $k_0$.
The condition means that $z$ is not contained in any proper
subvariety of $\PP$ defined over $k_0$ (in particular it is not in
the complement of $\PP')$.

\medskip
In the sequel, we will assume:

\medskip\noindent
(5.4.2)\,\,\, $D_0\,(=D_{z_0})$ is irreducible, and has only
moderate singularities (5.1).

\medskip\noindent
(5.4.3)\,\,\, The functor {\rm (1.1.1)} factors through $\MHM(X)$.

\medskip\noindent
{\bf 5.5.~Theorem.}
{\it With the notation and assumptions of {\rm (5.3)} and {\rm (5.4)},
assume further that if $j$ in {\rm (5.5.1)} is equal to $\dim X$,
then $k$ is a finitely generated subfield of $\C$ and the point
$z=(z_{1},\dots, z_{m})$ corresponding to $(D_{1},\dots, D_{m})$
belongs to $P'_0$ in {\rm (5.4)}.
If the image of {\rm (5.3.1)} is not one-dimensional, we assume also
that $(D_{1},\dots, D_{m})$ is $k_0$-generic for a subfield $k_0$ of
$k$ as in {\rm (5.4)}.
Let $X'=X\setminus\mcup_{0\le i\le m}\,D_i$.
Then the cycle map
$$
\CH^j(X',j)_{\Q}\to\Hom_{\M(k)}(\Q, H^j(X',\Q)(j))
\leqno(5.5.1)
$$
is surjective.
More precisely, the target of {\rm (5.5.1)} is generated by
$df_I/f_I$ for $I\subset\{1,\dots, r\}$ with $|I|=j$.
In particular, the target is zero if $j>r$.
}

\medskip\noindent
{\it Proof.} We proceed by increasing induction on $j$ and $m$.
We first consider the case where the image of (5.3.1) is
one-dimensional.
Then we may assume $\div f_i=a_0[D_i]-a_i[D_0]$.
Take $\xi$ from the target of (5.5.1), and consider its residue
along $D_{m}$
$$
\Res_{D_{m}}\xi\in\Hom_{\M(k)}(\Q, H^{j-1}(D_{m}
\setminus\mcup_{0\le i\le m-1}\,D_i,\Q)(j-1)).
$$
This is defined by using the connecting morphism of the localization
sequence, and at the level of logarithmic forms it is given by
residue.
Then the inductive hypothesis implies that $\Res_{D_{m}}\xi$ is a
linear combination of $df_I/f_I|_{D_{m}}$ for
$I\subset\{1,\dots, m-1\}$ such that $|I|=j-1$.
Note that the assertion for $j-1=0$ follows from (5.4.3) and the
assumption that $D_{m}$ is irreducible (by definition of $P'_0$
if $\dim X=1)$.
Indeed, (5.4.3) implies that the target of (5.5.1) is $\Q$ if
$j=0$ and $X'$ is irreducible.

Since $df_I/f_I|_{D_{m}}$ is the residue of
$a_0^{-1}df_{m}/f_{m}\wedge df_I/f_I$, we may assume that
$\Res_{D_{m}}\xi$ vanishes modifying $\xi$ by a linear combination
of products of integral logarithmic $1$-forms as above if necessary.
Then, using Lemma~(5.2), the assertion is reduced to the case where
$m$ is decreased by one, and we can proceed by induction.
In the case $m=0$, we have $\xi=0$ by (5.4.2) and (5.1.5).
So the assertion follows in the first case.

In the case where the image of (5.3.1) is not one-dimensional,
the argument is similar by using Lemma~(5.6) below.
We first show by induction on $j$ that the target of (5.5.1)
vanishes if (5.3.1) is injective.
For $j=1$, set $\tD=\coprod_iD_i$ and
$$
H^0(\tD,\Q)^0=\Ker(H^0(\tD,\Q)\to H^2(X,\Q)(1)).
$$
Here $H^0(\tD,\Q))=\mopl_i\,\Q[D_i]$ by Remark~(1.5)(v).
We have a short exact sequence
$$
0\to H^1(X,\Q)(1)\to H^1(X',\Q)(1)\to H^0(\tD,\Q)^0\to 0,
$$
so that the target of (5.5.1) is identified with the kernel of
$$
\Hom_{\M(k)}(\Q,H^0(\tD,\Q)^0)\to\Ext^1_{\M(k)}(\Q,H^1(X,\Q)),
$$
which coincides with the kernel of the Abel-Jacobi map for divisors.
So the assertion follows from the injectivity of (5.3.1).

If $j>1$, Lemma~(5.6) implies the injectivity of (5.3.1) with $X$
replaced by $D_i$ (and $k_0$ by the function field of the product
of $\PP_{L_j}$ for $j\ne i$ over $k_0$) for any $1\le i\le m$.
Hence the target of (5.5.1) for $D_i$ and $j-1$ vanishes by
inductive hypothesis.
Let $D=\mcup_{0\le i\le m}\,D_i$, $D'=D\setminus D_0$,
$Z=\Sing D'$, and $d=\dim X$.
Then the target of (5.5.1) vanishes by considering
$\Gr_{2j-2d}^{W}$ of the exact sequences
$$
\aligned
H_{2d-j}^{\BM}(X\setminus D_0,\Q)
&\to H_{2d-j}^{\BM}(X\setminus D,\Q)\to H_{2d-j-1}^{\BM}(D',\Q),
\\
H_{2d-j-1}^{\BM}(Z,\Q)
&\to H_{2d-j-1}^{\BM}(D',\Q)\to
H_{2d-j-1}^{\BM}(D'\setminus Z,\Q),
\endaligned
$$
since (5.1.5) holds for
$$
H_{2d-j}^{\BM}(X\setminus D_0,\Q)=H^j(X\setminus D_0,\Q)(d),
$$
and $\Gr_n^{W}H_i^{\BM}(Z,\Q)=0$ for $n>2\dim Z-2i$.
(We can verify the last vanishing by using the localization
sequence, because the smooth case is well-known [11].)

If (5.3.1) is not injective, we may assume that a multiple of
$[D_{m}]$ is rationally equivalent to a linear combination of
$[D_i]$ for $0\le i<m$, and we can apply the same argument as
in the first case by using Lemma~(5.6) and applying the inductive
hypothesis to $D_{m}$, where $k_0$ is replaced by the function
field of $\PP_{L_{m}}$ over $k_0$.
Thus the assertion is reduced to the case where $m$ is decreased
by one (using Lemma~(5.2)), and follows from the inductive
hypothesis.

Thus the proof of Theorem~(5.5) is reduced to the following:

\medskip\noindent
{\bf 5.6.~Lemma.} {\it Let $X$ be a smooth projective $k$-variety,
and $D_i$ be divisors on $X$ for $0\le i\le m$.
Assume $X$ and $D_i$ are defined over a subfield $k_0$ of $k$,
and $D_{m}$ is a $k_0$-generic hyperplane section of $X$.
Let $a_i\in\Z$, and assume
$\msum_{0\le i<m}\,a_i[D_i\cap D_{m}]=0$ in $\CH^1(D_{m})_{\Q}$.
Then $\msum_{0\le i<m}\,a_i[D_i]=0$ in $\CH^1(X)_{\Q}$.
}

\medskip\noindent
{\it Proof.} Let $X_{k_0}, D_{i,k_0}$ be $k_0$-varieties with
isomorphisms $X=X_{k_0}\otimes_{k_0}k$, $D_i=D_{i,k_0}\otimes_{k_0}k$.
Let $L$ be a very ample line bundle of $X_{k_0}$ such that $D_{m}$
corresponds to a $k$-valued generic point of $S_{k_0} :=\PP_{L}$ in
the notation of (5.4).
Let $Y$ denote the divisor on $X_{k_0}\times_{k_0}S_{k_0}$ whose
fiber over $z\in S_{k_0}(\ok_0)$ is the divisor corresponding to $z$.
By assumption, there exist a $k_0$-variety $S'_{k_0}$ and a dominant
morphism $\rho:S'_{k_0}\to S_{k_0}$ such that
$$
\msum_{0\le i<m}\,a_i[D'_{i,k}]=0\quad\text{in}\,\,\,\,
\CH^1(Y')_{\Q},
\leqno(5.6.1)
$$
where $Y'$ is the base change of $Y$ by $\rho$, and $D'_{i,k_0}$ is
the pull-back of $D_{i,k_0}$ to $Y'$ by the canonical morphism
$Y'\to X_{k_0}$.
Replacing $S_{k_0}$ with a locally closed subvariety if necessary,
we may assume that $\rho$ is\'etale, and then it is an open
embedding by using the pushforward under $\rho$.
Let $C_{k_0}$ be a generic line in $S_{k_0}$ which is not contained
in $S_{k_0}\setminus S'_{k_0}$.
Then the restriction of $Y\otimes_{k_0}k$ over
$C :=C_{k_0}\otimes_{k_0}k$ is a Lefschetz pencil $f :\tX\to C$, and
$\pi :\tX\to X$ is a blow-up along a smooth center $Z\subset X$ such
that $\codim Z=2$ and $D_i$ for $i<m$ intersects $Z$ properly
(in particular, $\pi^{*}D_i$ does not contain the exceptional divisor
of the blow-up $\pi)$.
Restricting (5.6.1) over a generic point of $C$, we get
$$
\msum_{0\le i<m}\,a_i[\pi^{*}D_i]=0\quad\text{in}\,\,\,\,
\CH^1(\tX)/f^{*}\CH^1(C),
$$
because the fibers of $f$ are irreducible.
Since $\CH^1(C)=\Z$, it implies
$$
\msum_{0\le i<m}\,a_i[\pi^{*}D_i]=c[\tX_{s}]\quad
\text{in}\,\,\,\,
\CH^1(\tX),
\leqno(5.6.2)
$$
where $\tX_{s}$ is the fiber of $f$ at a general $k$-valued point $s$
of $C$.
But this implies $c=0$ by applying $\pi^{*}\pi_*$ to (5.6.2), because
$\pi^{*}\pi_*\tX_{s}$ is the sum of $\tX_{s}$ and the exceptional
divisor of $\pi$ (and $\pi_*\pi^{*}=id)$.
So the assertion follows by applying $\pi_*$ to (5.6.2).
This completes the proofs of Lemma~(5.6) and Theorem~(5.5).

\medskip
By a similar argument, we can prove the Tate-type conjecture
corresponding to Theorem~(5.5).

\medskip\noindent
{\bf 5.7.~Theorem.} {\it With the notation and the assumptions of
Theorem~{\rm (5.5),} assume that $k$ is finitely generated.
Then the $l$-adic cycle map
$$
\CH^j(X',j)\otimes {\Q_{l}}\to
H^j(X'\otimes_{k}\ok,\Q_{l}(j))^{\Gal(\ok/k)}
\leqno(5.7.1)
$$
is surjective.
}

\medskip\noindent
{\it Proof.} It is well known that the $l$-adic Abel-Jacobi map
$$
\Pic_{Y}(k)^{0}\otimes_{\Z}\Q_{l}\to
H^1(\Gal(\ok/k),H^1(Y\otimes_{k}\ok,\Q_{l})(1))
$$
for divisors on a smooth proper $k$-variety $Y$ is expressed by
using the exact sequences associated with the Kummer sequence, and
is injective.
This can be used to show vanishing of the target of (5.7.1) when
(5.3.1) is injective and $j=1$.
The other arguments are similar to the proof of Theorem~(5.5).

\medskip\noindent
{\bf 5.8.~Theorem.} {\it With the notation and the assumptions of
{\rm (5.3)} and {\rm (5.4),} assume $(D_{1},\dots, D_{m})$ is
$k_0$-generic, and if $j$ in {\rm (5.8.1)} is equal to $\dim X$,
$D_I$ is not a rational curve for any $I\subset\{1,\dots, m\}$ such
that $|I|=\dim X-1$ \(including the case $\dim X=1$ and
$I=\emptyset)$.
Let $X'=X\setminus\mcup_{0\le i\le m}\,D_i$.
Then the cycle map
$$
\CH^j(X',j)_{\Q}\to\Hom_{\M(k)}(\Q, H^j(X',\Q)(j))
\leqno(5.8.1)
$$
is surjective.
More precisely, the target is generated by $df_I/f_I$ as in
Theorem~{\rm (5.5).} }

\medskip\noindent
{\it Proof.} We may assume $k=\C$ and $\M(k)=\MHS$ by
(5.4.3), because the assertion in the case $k=\C$ implies that
we have a desired higher cycle over a finite extension of $k$,
and its cycle class in the de Rham cohomology is invariant under
the action of the Galois group if $k$ is not algebraically closed.
Then the assertion follows from an argument similar to the proof of
Theorem~(5.5) by using the next Proposition (where only one variable
$z_i$ is free and the other $z_j$ are fixed) instead of Hilbert's
irreducibility theorem in the case $j=\dim X$.
When $j=1$, an element in the target of (5.8.1) is written as
$df/f$ with $f$ a rational function on $X$ if it comes from
integral cohomology, and $\div f=\Res\,df/f$ as well-known.
Then we can spread this function $f$ out so that we get a rational
function on $X\times S$ where $S$ is\'etale over $\PP'$.
Thus we get the spreading out of the given integral logarithmic
$1$-form as a horizontal family of integral logarithmic $1$-forms
over $S$.
The inductive argument is similar to the proof of Theorem~(5.5).

\medskip\noindent
{\bf 5.9.~Proposition.}
{\it With the notation and the assumptions of Theorem~{\rm (5.8)},
assume further $k=\C$, $m=1$ and $\dim X=1$ or $2$.
Let $I=\{1\}$ if $\dim X=1$, and $I=\{0,1\}$ if $\dim X=2$.
Let $\xi$ be a multivalued section of the local system
$\pi_*\Z_{D_I}$, where $\pi:D_I\to\PP'$ is as in {\rm (5.4).}
Let $\xi_{z}$ denote the $0$-cycle on $D_{I,z}$ defined by
the stalk of $\xi$ at $z\in\PP'(\C)$ \(which is also multivalued\).
Assume the image of $\xi_{z}$ in $\CH^1(X)$ is locally constant
for $z$ if $\dim X=1$, and that in $\CH^1(D_{1,z})$ vanishes
\(where $D_{1,z}=D_{z_{1}})$ if $\dim X=2$.
Then $\xi_{z}$ is a multiple of the canonical $0$-cycle $[D_{I,z}]$
\(i.e. the coefficient of $\xi_{z}$ at every point of
$D_{I,z}$ is same\) if $\dim X=1$, and $\xi=0$ if $\dim X=2$.
}

\medskip\noindent
{\it Proof.} We first show the case $\dim X=1$.
Let $n=\dim\PP$.
This coincides with the dimension of the projective space in which
$X$ is embedded by $L_{1}$.
We may assume $n\ge 2$ because the assertion is clear if $n=1$.
Since a hyperplane section is determined generically by $n$ points,
there exists an\'etale morphism of a non empty Zariski-open subset
$W$ of $X^n$ to $\PP'$.
Here we may assume that $W$ is stable by the action of the
permutation group on $X^n$.
Furthermore, there is a nonempty Zariski-open subset $U$ of $\PP'$
such that any ordered $n$ points of $D_{I,z}$ belongs to $W$ for
$z\in U(\C)$.
This is verified by using a finite \'etale morphism
$\sigma:S\to\PP'$ trivializing the monodromy group in
$\Aut(D_{I,z})$.
Indeed, for any subset $\Sigma\subset D_{I,z}$ with $|\Sigma|=n$,
let $B_{\Sigma}(\C)$ consist of points $s\in S$ such that the
parallel translate of $\Sigma$ over $s$ does not belong to $W$. Then
$U(\C)=\PP'(\C)\setminus\bigcup_{\Sigma}\sigma(B_{\Sigma}(\C))$.

Since $W$ is connected, this implies that the action of
$\pi_{1}(\PP'(\C),z)$ on $D_{I,z}$ is $n$-transitive for any
$z\in U(\C)$, and this holds for any $z\in\PP'(\C)$ since $U$ is
dense.
Then, for any two points of $D_{I,z}$, the image of
$\pi_{1}(\PP'(\C),z)$ in $\Aut(D_{I,z})$ contains a permutation of
$D_{I,z}$ which exchanges the given two points and keeps the other
points unchanged.
Indeed, this is easy for some two points using a Lefschetz pencil.
For any two points, we have a conjugate of it
using the $n$-transitivity since $n\ge 2$.

Assume the coefficients of $\xi_{z}$ at some two points of
$D_{I,z}$ are not same for $z\in\PP'(\C)$.
Then there exists $\rho\in\pi_{1}(\PP'(\C),z)$ such that the
coefficients of $\xi'_{z} :=\xi_{z}-\rho_*\xi_{z}$ are zero except
for two points of $D_{I,z}$ and $\xi'_{z}$ is a nonzero cycle of
degree $0$.
Furthermore the image of $\xi'_{z}$ in $\CH^1(X)$ is also locally
constant.
This implies that the map $X^{2}\to\Jac(X)$ defined by
$(x,y)\mapsto [x]-[y]$ is locally constant (using the above
morphism of $W$ to $\PP'$).
But this is clearly a contradiction (fixing $y$ for example).
So the assertion in the case $\dim X=1$ follows.

The argument is similar in the case $\dim X=2$.
Let $n=\dim\Gamma (D_0,L_{1}|_{D_0})-1$.
We may assume $n\ge 2$, because $D_0=\P^1$ and $D_{I,z}$ is one
point if $n=1$.
We can show the $n$-transitivity of the action of
$\pi_{1}(\PP'(\C),z)$ on $\Aut(D_{I,z})$ by the same argument as
above.
To construct a permutation of two points which keeps the other
points fixed, we consider a generic projection of $X$ to $\P^{2}$,
which is defined by choosing a generic three-dimensional subspace
of $\Gamma (X,L_{1})$.
Let $C_0$ be the image of $D_0$ in $\P^{2}$, and $C_{1}$ the
discriminant of $X\to\P^{2}$.
Then there is a hyperplane which is tangent to $C_0$ at a
sufficiently general smooth point (i.e. not at an inflection
point) of $C_0$ and the other intersections with $C=C_0\cup C_{1}$
are transversal.
Take the pull-back $D_{1}$ of the hyperplane to $X$, and consider a
pencil containing $D_{1}$.
(Note that $D_{1}$ is smooth because the above hyperplane intersects
$C_{1}$ transversely.)

Let $\Delta$ be a sufficiently small open disk in the base space
$\P^1$ of the pencil (which is viewed as a subset of $\PP'(\C)$)
such that the fiber at the origin is $D_{1}$.
There is a connected component $\Delta'$ of
$\mcup_{z\in\Delta}\,D_0\cap D_{1,z}$ which is a ramified covering
of degree $2$ over $\Delta$, and the other connected components are
biholomorphic to $\Delta$.
If the coefficients of $\xi_{z}$ at some two points of
$D_0\cap D_{1,z}$ are not same, we may assume that these two points
are $\Delta'\cap D_{1,z}$ by the $n$-transitivity.
Then, using the local monodromy around $0\in\Delta$, we get a
nonzero cycle $\xi'_{z}$ which is supported on $\Delta'\cap D_{1,z}$,
and has degree $0$, and whose image in $\Jac(D_{1,z})$ vanishes.
Taking the base change by $\Delta'\to\Delta$, we get a family of
smooth proper curves $D_{1,z'}\,(z'\in\Delta')$ which has two
sections $s$ and $s'$ such that
$\Im\,s\cup\Im\,s'=\Delta'\times_{\Delta}\Delta'$.
This gives a univalent family of nonzero cycles $\xi'_{z'}$
supported on $\Delta'\times_{\Delta}\{z'\}$ for $z'\in\Delta'$
such that the image of $\xi'_{z'}$ in $\Jac(D_{1,z'})$ vanishes.
This induces a contradiction by considering the embedding
$D_{1,z'}\to\Jac(D_{1,z'})$ determined by $s_{z'}$.
So we get $\xi=0$ because $\xi_{z}$ has degree zero.
This finishes the proofs of Proposition~(5.9) and Theorem~(5.8).

\medskip
By an argument similar to the proofs of Theorem~(5.8) and
Proposition~(5.9), we can show the following (which is compatible
with Voisin's conjecture [35], see (4.4)):

\medskip\noindent
{\bf 5.10.~Proposition.}
{\it With the notation and the assumptions of Theorem~{\rm (5.5)} or
{\rm (5.8),} let $D=\mcup_{0\le i\le m}\,D_i$, and $\tD$ be the
normalization of $D$.
Then the morphism
$$
\CH^{2}(X',2)_{\Q}\oplus\CH^1(\tD,1)_{\Q}\to\CH^1(D,1)_{\Q}
\leqno(5.10.1)
$$
is surjective.
More precisely, we have the surjection with $\CH^{2}(X',2)$
replaced by the second Milnor $K$-group of
$\Gamma (X',{\O}_{X'}^{*})$, and the morphism to the target is
given by the tame symbol.
In particular, there is no nontrivial indecomposable higher cycle
in $\CH_{\ind}^{2}(X,1)_{\Q}$ which is supported on $D$.
}

\medskip\noindent
{\it Proof.} This follows from Theorem~(5.5) or (5.8) by increasing
induction on $m$.
Indeed, let $D'_i=D_i\setminus\mcup_{j\ne i} D_j$.
The kernel of the cycle map (5.5.1) or (5.8.1) for $D'_{m}$ with
$j=1$ comes from $\CH^1(D_{m},1)\,(=\Gamma(D_{m},\O_{D_{m}}^{*}))$,
because $\div f=0$ on $D_{m}$ if
$f\in\Gamma(D'_{m},\O_{D'_{m}}^{*})$ belongs to the kernel.
The assertion is then reduced to the case with $m$ decreased by one
(using Theorem~(5.5) or (5.8)), because the logarithmic differential
of the tame symbol $\{f,g\}$ is given by the residue of
$df/f\wedge dg/g$ up to a sign.
In the case $m=0$, we have $\CH^1(D_0,1)=\CH^1(\tD_0,1)$ because
$D_0$ is analytic-locally irreducible.
So the assertion follows.

\medskip\noindent
{\bf 5.11.~Theorem.} {\it Let $X$ be a smooth projective
$k$-variety with an ample line bundle $L$, and $X'$ be the
complement of a union of smooth hypersurfaces $D_0,\dots, D_{m}$
with respect to $L$ such that $D :=\mcup_{0\le i\le m} D_i$ is a
divisor with normal crossings.
Let $p, n$ be positive integers such that $p>n$ and $2p-n<\dim X$.
Assume that the $\M$-Hodge conjecture is true for codimension
$p-n$ cycles on $X$.
Then the cycle map
$$
\CH^p(X',n)_{\Q}\to
\Hom_{\M(k)}(\Q,H^{2p-n}(X',\Q)(p))
\leqno(5.11.1)
$$
is surjective.
More precisely, the target is generated by
$df_I/f_I\wedge cl(\zeta)$ in the notation of {\rm (5.3)} where
$I$ is a subset of $\{1,\dots, r\}$ with $|I|=n$, and
$\zeta\in\CH^{p-n}(X)$.
}

\medskip\noindent
{\it Proof.} This follows from an argument similar to the proofs of
Theorems~(5.5) and (5.8) by increasing induction on $m$ and $j$.
Indeed, the condition $2p-n<\dim X$ implies that $2p-2n<\dim D_I$,
and the restriction morphism induces an isomorphism
$$
H^{2p-2n}(X,\Q)\to H^{2p-2n}(D_I,\Q)
$$
by the weak Lefschetz theorem.
Furthermore, the condition implies in the case $m=0$ that
$H^{2p-n}(X',\Q)$ is pure of weight $2p-n$ by an argument similar
to the proof of Lemma~(5.2) (using the hard Lefschetz theorem).
This finishes the proof of Theorem~(5.11).

\medskip\noindent
{\bf 5.12.~Remarks.} (i) The last assertion of (5.12) does not hold
unless the $D_i$ are hypersurface sections of the same ample line
bundle.
For example consider $X=\P^1\times\P^{2}$ with $D_i$ a general
hyperplane section of $\O_{\P^1}(a_i)\otimes\O_{\P^{2}}(b_i)$
for $i=0, 1$, where $a_i, b_i$ are positive integers such that
$a_0b_{1}-a_{1}b_0\ne 0$.
Then (5.3.1) is injective, i.e. $r=0$.
But the target of (5.1.1) does not vanish for $p=2, n=1$.
Indeed, $\dim H^{4}(X,\Q)=2$, and the Gysin morphism
$H^{2}(D_i,\Q)\to H^{4}(X,\Q)(1)$ is an isomorphism, because the
restriction morphism $H^{2}(X,\Q)\to H^{2}(D_i,\Q)$ and its
composition with the Gysin morphism are isomorphisms by the hard
and weak Lefschetz theorems.
Using the weight spectral sequence, this implies that
$\Gr_{4}^{W} H^{3}(X',\Q)=\Q(-2)$, and the assertion follows
because $\Gr_{3}^{W} H^{3}(X',\Q)=0$.

\medskip
(ii) The converse of Theorem~(4.1) is not true in general even if
(4.1.1) is restricted to higher cycles supported on the complement
of $U$ for which the surjectivity of (4.1.2) holds.
Indeed, let $X$ be a smooth complex projective variety such that
$\Gamma(X,\Omega_{X}^1)\ne 0$.
Take very ample line bundles $L_0, L_{1}$ such that
$L_0\otimes L_{1}^{\vee}$ is a non-torsion point of the Picard
variety of $X$.
Let $D_0, D_{1}$ be general hyperplane sections of $L_0, L_{1}$
which intersect transversely.
Let $X'=X\setminus (D_0\cup D_{1})$.
Then, by an argument similar to the proof of Proposition~(5.9),
we get
$$
\Hom_{\MHS}(\Q,H^{2}(X',\Q)(2))=0 .
$$
In particular, (4.1.2) is surjective for $n=1, p=2$, and the
surjectivity of (4.1.3) is easy, see [20].
Consider now a decomposable higher cycle
$\zeta :=([D_0]-[D_{1}])\otimes\alpha\in\CH^{2}(X,1)_{\Q}$.
It is nonzero if $\alpha\in\C$ is not algebraic over a subfield
$k$ on which $X$ and $D_i$ are defined (see [32]).
But $cl(\zeta)$ always vanishes because $[D_0]-[D_{1}]$ is
homologically equivalent to zero.

\end{document}